\documentclass[12pt]{amsart}
\usepackage{amsmath, amsthm, amscd, amsfonts}

\newtheorem{theorem}{Theorem}[section]
\newtheorem{lemma}[theorem]{Lemma}
\newtheorem{proposition}[theorem]{Proposition}
\newtheorem{corollary}[theorem]{Corollary}
\theoremstyle{definition}
\newtheorem{definition}[theorem]{Definition}
\newtheorem{example}[theorem]{Example}

\theoremstyle{remark}
\newtheorem{remark}[theorem]{Remark}
\numberwithin{equation}{section}
\begin{document}
\title {$\sigma$-Derivations in Banach Algebras}
\author [M. Mirzavaziri, M. S. Moslehian]{Madjid Mirzavaziri$^1$ and Mohammad Sal Moslehian$^2$}
\address{$^{1}$ Department of Mathematics, Ferdowsi University, P. O. Box 1159, Mashhad 91775,
Iran; \newline Centre of Excellency in Analysis on Algebraic
Structures (CEAAS).} \email{mirzavaziri@math.um.ac.ir}
\address{$^{2}$ Department of Mathematics, Ferdowsi University, P. O. Box 1159, Mashhad 91775, Iran;
\newline Centre of Excellence in Analysis
on Algebraic Structures (CEAAS).}
\email{moslehian@ferdowsi.um.ac.ir} \subjclass[2000]{Primary
47D06; Secondary 47B47, 16W25}
\keywords{(inner)
$\sigma$-derivation; (inner) $\sigma$-endomorphism; one-parameter
group of $\sigma$-endomorphisms; generalized Leibniz rule.}

\begin{abstract} Introducing the notions of (inner) $\sigma$-derivation,
(inner) $\sigma$-endomorphism and one-parameter group of
$\sigma$-endomorphisms ($\sigma$-dynamics) on a Banach algebra,
we correspond to each $\sigma$-dynamics a $\sigma$-derivation
named as its $\sigma$-infinitesimal generator. We show that the
$\sigma$-infinitesimal generator of a $\sigma$-dynamics of inner
$\sigma$-endomorphisms is an inner $\sigma$-derivation and deal
with the converse. We also establish a nice generalized Leibniz
formula and extend the Kleinenckr--Sirokov theorem for
$\sigma$-derivations under certain conditions.
\end{abstract}
\maketitle{}

\section{Introduction}
\noindent Let ${\mathcal A}$ be a Banach algebra. Recall that a
derivation $d$ defined on a (dense) subalgebra ${\mathcal D}$ of
${\mathcal A}$ is a linear mapping satisfying $d(ab)=d(a)b+ad(b),
a,b\in {\mathcal D}$. A derivation $d$ is said to be inner if there
exists an element $u\in {\mathcal A}$ such that $d(a)=ua-au$ for all
$a\in {\mathcal A}$. There are nonzero derivations defined on a
commutative algebra among which we may consider the ordinary
derivative $d/dt:C^1([0,1])\to C([0,1])$ where $C^1([0,1])$ is the
algebra of all continuously differentiable functions on $[0,1]$.
This example gives an idea to define a derivation on a dense
subalgebra of a given algebra ${\mathcal A}$.

Derivations play essential role in some important branches of
mathematics and physics such as dynamical systems. The general
theory of dynamical systems is the paradigm for modeling and
studying phenomena that undergo spatial and temporal evolution. The
application of dynamical systems has nowadays spread to a wide
spectrum of disciplines including physics, chemistry, biochemistry,
biology, economy and even sociology. In particular, the theory of
dynamical systems concerns the theory of derivations in Banach
algebras and is motivated by questions in quantum physics and
statistical mechanics, cf. \cite{SAK}.

It is known that the relation $TS-ST=I$ is impossible for bounded
operators $T$ and $S$ on Banach spaces; cf. \cite{SAK} and
references therein. In fact the study of this relation as a special
case of $T\sigma(S)-\sigma(S)T=R$, where $\sigma$ is a linear
mapping, leads the theory of derivations to be extensively
developed.

The above considerations motivate us to generalize the notion of
derivation as follow. Let ${\mathcal D}$ be a subalgebra of a Banach
algebra ${\mathcal A}$ and let $\sigma,d:{\mathcal D}\to {\mathcal
A}$ be linear mappings. If $d(ab)=d(a)\sigma(b)+\sigma(a)d(b)$ for
all $a,b\in {\mathcal D}$ then we say $d$ is a $\sigma$-derivation
(see \cite{A-M-N, A-R, B-M, H-K, KOM, MOS1, MOS2} and references
therein). There are some interesting questions in this area of
research, e.g. one may ask `What are the $\sigma$-derivations of the
compact operators acting on a separable Hilbert space?' The paper
\cite{BAT} can be a starting point for answering this question. Note
that if $\sigma$ is the identity map then $d$ (and every so-called
inner $\sigma$-derivation) is indeed a derivation (inner derivation,
respectively) in the usual sense.

In this paper we introduce and study (inner) $\sigma$-derivations,
(inner) $\sigma$-endomorphisms and one-parameter group of
$\sigma$-endomorphisms ($\sigma$-dynamics). The importance of our
approach is that $\sigma$ is a linear mapping, not necessary an
algebra endomorphism. It is shown that the $\sigma$-infinitesimal
generator of a $\sigma$-dynamics of inner $\sigma$-endomorphisms is
an inner $\sigma$-derivation and the converse is true under some
conditions. We give a formula for computation $d^n(ab)$ where $d$ is
a $\sigma$-derivation that is interesting in its own right. We also
generalize two known theorems in the context of Banach algebras,
namely the Wielandt--Wintner theorem and Kleinecke--Shirokov
theorem.

This paper is self-contained. The reader, however, is referred to
\cite{DAL} for details on Banach algebras and to \cite{B-R1, B-R2,
SAK} for more information on dynamical systems.

\section{$\sigma$-dynamics}
\noindent Throughout the paper ${\mathcal A}$ denotes a Banach
algebra, $\iota$ is the identity operator on ${\mathcal A}$,
${\mathcal D}$ denotes a subalgebra of ${\mathcal A}$, and $\sigma,
d:{\mathcal D}\to {\mathcal A}$ are linear mappings.

\begin{definition} $d$ is called a $\sigma$-derivation if $d(ab)=d(a)\sigma(b)+\sigma(a)d(b)$
for all $a,b\in {\mathcal D}$.\end{definition}

\begin{example} Let $\sigma$ be an arbitrary linear mapping on ${\mathcal D}$ and suppose that
$u$ is an element of ${\mathcal A}$ satisfying
$u(\sigma(ab)-\sigma(a)\sigma(b))=(\sigma(ab)-\sigma(a)\sigma(b))u$
for all $a,b\in {\mathcal D}$. Then the mapping $d:{\mathcal D}\to
{\mathcal A}$ defined by $d(a)=u\sigma(a)-\sigma(a)u$ is a
$\sigma$-derivation.\end{example}

The above $\sigma$-derivation is called inner. Note that if
$\sigma$ is a endomorphism then $u$ can be any arbitrary element
of ${\mathcal A}$.

\begin{example} Let $d$ be an endomorphism on ${\mathcal A}$. Then $d$ is a
$\frac{d}{2}$-derivation.\end{example}

\begin{example} Let $d,\sigma :C([0,1])\to C([0,1])$ be defined by $\sigma(f)=\frac{f}{2}$
and $d(f)=fh_0$, respectively. Here $h_0$ is an arbitrary fixed
element in $C([0,1])$. Then easy observations show that
$\sigma(1)\neq 1$, $d(1)\neq 0$, the linear mapping $\sigma$ is
not endomorphism, and $d$ is a $\sigma$-derivation.\end{example}

\begin{example} Suppose that $d,\sigma:C([0,1])\to C([0,1])$ are defined by
$\sigma(f)(t)=\frac{1}{2}f(2t)$ and $d(f)(t)=f(2t)h_0(t)$,
respectively, where $h_0$ is an arbitrary fixed element of
$C([0,1])$. Then
\begin{eqnarray*}
d(fg)(t)&=&f(2t)g(2t)h_0(t)\\
&=&(f(2t)h_0(t))(\frac{1}{2}g(2t))+(\frac{1}{2}f(2t))(g(2t)h_0(t))\\
&=&d(f)(t)\sigma(g)(t)+\sigma(f)(t)d(g)(t)
\end{eqnarray*}
It follows that $d$ is a $\sigma$-derivation and no scalar
multiple of $\sigma$ is a endomorphism.
\end{example}

\begin{definition} A linear mapping $\alpha:{\mathcal A}\to {\mathcal A}$ is called
$\sigma$-endomorphism if
$(\alpha+\sigma-\iota)(ab)-(\alpha+\sigma-\iota)(a)(\alpha+\sigma-\iota)(b)=\sigma(ab)-\sigma(a)\sigma(b)$
for all $a,b\in {\mathcal A}$.\end{definition}

Note that if $\sigma=\iota$ then a $\sigma$-endomorphism is nothing
more than an endomorphism on ${\mathcal A}$ in the usual sense.

\begin{lemma}
Let $\alpha$ be a linear mapping on ${\mathcal A}$. Then $\alpha$
is a $\sigma$-endomorphism if and only if
$$\alpha(ab)-\alpha(a)\alpha(b)=(\alpha(a)-a)(\sigma(b)-b)+(\sigma(a)-a)(\alpha(b)-b)$$
\end{lemma}
\begin{proof} Straightforward.
\end{proof}

\begin{definition} A mapping $t\in {\mathbb R} \mapsto \alpha_t\in
B({\mathcal A})$ denoted by $\{\alpha_t\}_{t\in{\mathbb R}}$ is a
one-parameter group of bounded operators on ${\mathcal A}$ if it
satisfies the following conditions:

(i) $\alpha_t\alpha_s=\alpha_{t+s}$, for all $t,s\in{\mathbb R}$;

(ii) $\alpha_0=\iota$.

In the case that $\alpha_t$'s are bounded $\sigma$-endomorphisms,
$\{\alpha_t\}_{t\in{\mathbb R}}$ is called a one-parameter group
of $\sigma$-endomorphisms on ${\mathcal A}$. It is said to be
uniformly continuous if the map $t\mapsto\alpha_t$ is continuous
in the uniform topology, i.e. $\|\alpha_t-\iota\|\to 0$ as $t\to
0$. In this case $\{{\mathcal A},\alpha\}$ is called a
$\sigma$-dynamics.\end{definition}

Let $\{{\mathcal A},\alpha\}$ be a $\sigma$-dynamics. Then for
each $a\in {\mathcal A}$, if the limit of
$t^{-1}(\alpha_t(a)-\iota(a))$ as $t$ tends to $0$ exists, we can
define $d(a)$ to be this limit. This provides a mapping
$d:{\mathcal D}\to {\mathcal A}$, where ${\mathcal D}$ is the set
of all elements $a$ in ${\mathcal A}$ for which the limit exists.
$d$ is called the $\sigma$-infinitesimal generator of the
$\sigma$-dynamics $\{\alpha_t\}_{t\in{\mathbb R}}$.

\begin{proposition} Let $\{{\mathcal A},\alpha\}$ be a $\sigma$-dynamics.
Then $d=\displaystyle{\lim_{t\to 0}}t^{-1}(\alpha_t(a)-a)$ is an
everywhere defined $\sigma$-derivation.\end{proposition}
\begin{proof} We have
\begin{eqnarray*}
d(ab)&=&\lim_{t\to 0}t^{-1}(\alpha_t(ab)-\iota(ab))\\
&=&\lim_{t\to 0}t^{-1}((\alpha_t+\sigma-\iota)(ab)-\sigma(ab))\\
&=&\lim_{t\to
0}t^{-1}((\alpha_t+\sigma-\iota)(a)(\alpha_t+\sigma-\iota)(b)-\sigma(a)\sigma(b))\\
&=&\lim_{t\to 0}(t^{-1}((\alpha_t+\sigma-\iota)(a)-\sigma(b))+ (\alpha_t+\sigma-\iota)(a)t^{-1}((\alpha_t+\sigma-\iota)(b)
-\sigma(b)))\\
&=&d(a)\sigma(b)+\sigma(a)d(b).
\end{eqnarray*}
It follows from Proposition 3.1.1 of \cite{B-R1} that $d$ is
everywhere defined.
\end{proof}

\begin{definition} A linear mapping $\alpha$ is called an inner $\sigma$-endomorphism if there is an
element $u\in {\mathcal A}$ such that
$(\alpha+\sigma-\iota)(a)=e^{u}\sigma(a)e^{-u}, a\in {\mathcal
A}$ and
$u(\sigma(ab)-\sigma(a)\sigma(b))=(\sigma(ab)-\sigma(a)\sigma(b))u,
a,b\in {\mathcal A}$.\end{definition}

\begin{lemma} Each inner $\sigma$-endomorphism is indeed a $\sigma$-endomorphism.
\end{lemma}
\begin{proof}
\begin{eqnarray*}
(\alpha+\sigma-\iota)(ab)-(\alpha+\sigma-\iota)(a)(\alpha+\sigma-\iota)(b)&=&e^{u}(\sigma(ab)-\sigma(a)\sigma(b))e^{-u}\\
&=&e^{u}e^{-u}(\sigma(ab)-\sigma(a)\sigma(b))\\
&=&\sigma(ab)-\sigma(a)\sigma(b).
\end{eqnarray*}
\end{proof}

\begin{theorem} Let $\{\alpha_t\}_{t\in{\mathbb R}}$ be a one-parameter group of inner $\sigma$-endomorphisms.
Then the $\sigma$-infinitesimal generator $d$ of
$\sigma$-dynamics $({\mathcal A},\alpha)$ is an inner
$\sigma$-derivation.\end{theorem}
\begin{proof}
\begin{eqnarray*}
\lim_{t\to 0}t^{-1}(\alpha_t(a)-a)&=&\lim_{t\to 0}t^{-1}((\alpha_t+\sigma-\iota)(a)-\sigma(a))\\
&=&\lim_{t\to 0}t^{-1}(e^{tu}\sigma(a)e^{-tu}-\sigma(a))\\
&=&\lim_{t\to 0}(ue^{tu}\sigma(a)e^{-tu}-e^{tu}\sigma(a)ue^{-tu})\\
&=&u\sigma(a)+\sigma(a)u
\end{eqnarray*}
Note that we use L'Hospital's rule to get the third equality.
\end{proof}

\begin{lemma} Let $d:{\mathcal A}\to {\mathcal A}$ be the inner $\sigma$-derivation
$d(a)=u\sigma(a)-\sigma(a)u$. If $\sigma^2=\sigma$ and
$\sigma(au)=\sigma(a)u, \sigma(ua)=u\sigma(a)$ for all $a\in
{\mathcal A}$, then
\[\sum_{k=0}^r(-1)^k{r\choose k}u^k\sigma(a)u^{r-k}=(-1)^rd^r(a)\;\;\;\;(*)\] for all $a\in
{\mathcal A}, 0\leq k \leq r, r\geq 1$.
\end{lemma}
\begin{proof} We use induction on $r$. For $r=1$ there is nothing to do. Assume that $(*)$ holds for $r$. We have
\begin{eqnarray*}
(-1)^{r+1}d^{r+1}(a)&=&(-1)^{r+1}d(d^r(a))\\
&=&-u\sigma((-1)^{r}d^r(a))+\sigma((-1)^rd^r(a))u\\
&=&-\sum_{k=0}^r(-1)^k{r\choose k}u^{k+1}\sigma(a)u^{r-k}+\sum_{k=0}^r(-1)^k{r\choose k}u^k\sigma(a)u^{r-k+1}\\
&=&-\sum_{k=0}^r(-1)^k{r\choose k}u^{k+1}\sigma(a)u^{r-k}-\sum_{k=-1}^{r-1}(-1)^{k+1}{r\choose k+1}u^{k+1}\sigma(a)u^{r-k}\\
&=&(-1)^{r+1}u^{r+1}\sigma(a)+\sigma(a)u^{r+1}-\sum_{k=0}^{r-1}(-1)^k({r\choose k}+{r\choose k+1})u^{k+1}\sigma(a)u^{r-k}\\
&=&(-1)^{r+1}u^{r+1}\sigma(a)+\sigma(a)u^{r+1}-\sum_{k=0}^{r-1}(-1)^k{r+1 \choose k+1}u^{k+1}\sigma(a)u^{r-k}\\
&=&(-1)^{r+1}u^{r+1}\sigma(a)+\sigma(a)u^{r+1}+\sum_{k=1}^r(-1)^k{r+1 \choose k}u^k\sigma(a)u^{r+1-k}\\
&=&\sum_{k=0}^{r+1}(-1)^k{r+1 \choose k}u^k\sigma(a)u^{r+1-k}.
\end{eqnarray*}
\end{proof}

\begin{theorem} Let $d:{\mathcal A}\to {\mathcal A}$ be the inner $\sigma$-derivation
$d(a)=u\sigma(a)-\sigma(a)u$. If $\sigma^2=\sigma$ and $\sigma(au)=\sigma(a)u, \sigma(ua)=u\sigma(a)$ for
all $a\in {\mathcal A}$, then there exists a one-parameter group of operators $\{\alpha_t\}_{t\in{\mathbb R}}$ such
that $d$ is its $\sigma$-infinitesimal generator and $\alpha_t-\sigma+\iota$ is inner $\sigma$-homomorphism for
all $t\in{\mathbb R}$.\end{theorem}
\begin{proof} Put $\alpha_t(a)=\displaystyle{\sum_{n=0}^\infty}\frac{t^nd^n(a)}{n!}$. Using Lemma 2.13 we have
\begin{eqnarray*}
e^{tu}\sigma(a)e^{-tu}&=&(\sum_{n=0}^\infty\frac{t^nu^n}{n!})\sigma(a)(\sum_{m=0}^\infty\frac{(-t)^mu^m}{m!})\\
&=&\sum_{n=0}^\infty \sum_{m=0}^\infty\frac{t^n(-t)^m}{n!m!}u^n\sigma(a)u^m\\
&=&\sum_{r=0}^\infty \sum_{k=0}^r\frac{t^k(-t)^{r-k}}{k!(r-k)!}u^k\sigma(a)u^{r-k}\\
&=&\sum_{r=0}^\infty\frac{(-t)^r}{r!} \sum_{k=0}^r(-1)^k{r\choose k}u^k\sigma(a)u^{r-k}\\
&=&\sum_{r=0}^\infty\frac{t^r}{r!}d^r(a)\\
&=&\alpha_t(a).
\end{eqnarray*}
Since
$((\alpha_t-\sigma+\iota)+\sigma-\iota)(a)=\alpha_t(a)=e^{tu}\sigma(a)e^{-tu}$,
we deduce that $\alpha_t-\sigma+\iota$ is an inner
$\sigma$-endomorphism. Obviously $\alpha_t\alpha_s=\alpha_{t+s}$
and $\alpha_0=\iota$. In addition $\displaystyle{\lim_{t\to
0}}\frac{\alpha_t(a)-a}{t}=\displaystyle{\lim_{t\to
0}}\displaystyle{\sum_{n=1}^\infty}\frac{t^nd^n(a)}{n!}=d(a)$.
\end{proof}

\begin{definition} Let $d$ be a $\sigma$-derivation. We say $d$ multiplizes $\sigma$ if
$\sigma(ab)-\sigma(a)\sigma(b)\subseteq \ker(d)$. In this case $d$
is called a multiplizing $\sigma$-derivation.\end{definition}

\begin{example}Each inner $\sigma$-derivation $d$ is multiplizing. Let
$d(a)=u\sigma(a)-\sigma(a)u$ for some $u\in {\mathcal A}$. Then
we have $ d(ab)=d(a)\sigma(b)+\sigma(a)d(b)$, and so
$u\sigma(ab)-\sigma(ab)u=(u\sigma(a)-\sigma(a))u+\sigma(a)(u\sigma(b)-\sigma(b)u)$,
which implies
$u(\sigma(ab)-\sigma(a)\sigma(b))-(\sigma(ab)-\sigma(a)\sigma(b))u=0$
or $d(\sigma(ab)-\sigma(a)\sigma(b))=0$. Thus $d$ is a
multiplizing $\sigma$-derivation.\end{example}

\begin{proposition} Let ${\mathcal A}$ be an algebra with no zero divisor. Then $d$ is a
multiplizing $\sigma$-derivation if and only if
$\sigma(b\sigma(ab))=\sigma(b)\sigma^2(ab)$ for all $a,b\in
{\mathcal A}$.\end{proposition}
\begin{proof} For each $a,b,c\in{\mathcal D}$ we have

\begin{eqnarray*} d(abc)&=&d(ab)\sigma(c)+\sigma(ab)d(c)\\
&=&(d(a)\sigma(b)+\sigma(a)d(b))\sigma(c)+\sigma(ab)d(c)\\
&=&d(a)\sigma(b)\sigma(c)+\sigma(a)d(b)\sigma(c)+\sigma(ab)d(c).\end{eqnarray*}
On the other hand,
\begin{eqnarray*} d(abc)&=&d(a)\sigma(bc)+\sigma(a)d(bc)\\
&=&d(a)\sigma(bc)+\sigma(a)(d(b)\sigma(c)+\sigma(b)d(c))\\
&=&d(a)\sigma(bc)+\sigma(a)d(b)\sigma(c)+\sigma(a)\sigma(b)d(c).\end{eqnarray*}
Therefore
\[
d(a)(\sigma(bc)-\sigma(b)\sigma(c))=(\sigma(ab)-\sigma(a)\sigma(b))d(c)\]
for each $a,b,c\in{\mathcal D}$. Putting
$c=\sigma(ab)-\sigma(a)\sigma(b)$ we have
$\sigma(ab)-\sigma(a)\sigma(b)\subseteq \ker(d)$ if and only if
$\sigma(b(\sigma(ab)-\sigma(a)\sigma(b))-\sigma(b)\sigma(\sigma(ab)-\sigma(a)\sigma(b))=0$,
which implies the result.
\end{proof}

\section{Generalized Leibniz rule}
\noindent In the rest of the paper we need a family of mappings
$\{\varphi_{n,k}\}_{n\in N,0\leq k\leq 2^n-1}$ to simplify the
notations. We introduce these mappings by presenting the natural
numbers in base $2$.

Let $n$ be a natural number and $0\leq k\leq 2^n-1$. Note that
$2^n-1=(\underbrace{1\ldots 1}_{n~times})_2$ and each $0\leq
k\leq 2^n-1$ has at most $n$ digits in base $2$. Now assume that
$\varphi_{n,k}$ is the mapping derived from writing $k$ in base
$2$ with exactly $n$ digits and put $d$ for $1$'s and $\sigma$
for $0$'s.

To illustrate the mappings $\varphi_{n,k}$'s, let us give an
example. Let $n=5$ and $k=11$. Then we can write $k=(01011)_2$
and so $\varphi_{5,11}=\sigma d\sigma d d=\sigma d\sigma d^2$.
\begin{lemma} Let $n$ be a natural number and $0\leq k\leq 2^n-1$. Then

(i) $d\varphi_{n,k}=\varphi_{n+1,2^n+k}$;

(ii) $d\varphi_{n,2^n-1-k}=\varphi_{n+1,2^{n+1}-1-k}$;

(iii) $\sigma\varphi_{n,k}=\varphi_{n+1,k}$;

(iv)
$\sigma\varphi_{n,2^n-1-k}=\varphi_{n+1,2^{n+1}-1-(2^n+k)}$.\end{lemma}
\begin{proof} Assume that $k=(c_n\ldots c_1)_2$.

(i) $d\varphi_{n,k}=\varphi_{n+1,(1c_n\ldots
c_1)_2}=\varphi_{n+1,2^n+k}$.

(ii) $2^n-1-k=(\bar{c_n}\ldots\bar{c_1})_2$, where
$\bar{c_i}+c_i=1$, since $(2^n-1-k)+k=2^n=(1\ldots 1)_2$. Thus we
infer that
$d\varphi_{n,2^n-1-k}=d\varphi_{n,(\bar{c_n}\ldots\bar{c_1})_2}=
\varphi_{n+1,(1\bar{c_n}\ldots\bar{c_1})_2}=\varphi_{n+1,2^n+(2^n-1-k)}
=\varphi_{n+1,2^{n+1}-1-k}$.

(iii) $\sigma\varphi_{n,k}=\varphi_{n+1,(0c_n\ldots
c_1)_2}=\varphi_{n+1,k}$.

(iv)
$\sigma\varphi_{n,2^n-1-k}=\varphi_{n+1,(0\bar{c_n}\ldots\bar{c_1})_2}=
\varphi_{n+1,2^n-1-k}=\varphi_{n+1,2^{n+1}-1-(2^n+k)}$.
\end{proof}

\begin{theorem} For each $a,b\in{\mathcal D}$,
\[ d^n(ab)=\sum_{k=0}^{2^n-1}\varphi_{n,k}(a)\varphi_{n,2^n-1-k}(b)\hspace{2cm}(**)\]
\end{theorem}
\begin{proof} We prove the assertion by induction on $n$.
For $n=1$ we have
\[
d(ab)=d(a)\sigma(b)+\sigma(a)d(b)=
\varphi_{1,1}(a)\varphi_{1,0}(b)+\varphi_{1,0}(a)\varphi_{1,1}(b).\]
Now suppose $(**)$ is true for $n$. By Lemma 3.1 we obtain
\begin{eqnarray*}
d^{n+1}(ab)&=&d(d^{n}(ab))=d(\sum_{k=0}^{2^n-1}\varphi_{n,k}(a)\varphi_{n,2^n-1-k}(b))\\
&=&\sum_{k=0}^{2^n-1}d(\varphi_{n,k}(a)\varphi_{n,2^n-1-k}(b))\\
&=&\sum_{k=0}^{2^n-1}d(\varphi_{n,k}(a))\sigma(\varphi_{n,2^n-1-k}(b))+
\sigma(\varphi_{n,k}(a))d(\varphi_{n,2^n-1-k}(b))\\
&=&\sum_{k=0}^{2^n-1}\varphi_{n+1,2^n+k}(a)\varphi_{n+1,2^{n+1}-1-(2^n+k)}(b)\\
&&+\sum_{k=0}^{2^n-1}\varphi_{n+1,k}(a)\varphi_{n+1,2^{n+1}-1-k}(b)\\
&=&\sum_{l=2^n}^{2^{n+1}-1}\varphi_{n+1,l}(a)\varphi_{n+1,2^{n+1}-1-l}(b)\\
&&+\sum_{l=0}^{2^n-1}\varphi_{n+1,l}(a)\varphi_{n+1,2^{n+1}-1-l}(b)\\
&=&\sum_{l=0}^{2^{n+1}-1}\varphi_{n+1,l}(a)\varphi_{n+1,2^{n+1}-1-l}(b).
\end{eqnarray*}
\end{proof}

\begin{example} As an illustration, consider $n=3$. We have
\begin{eqnarray*}d^3(ab)&=&\varphi_{3,0}(a)\varphi_{3,7}(b)+\varphi_{3,1}(a)\varphi_{3,6}(b)\\
&&+\varphi_{3,2}(a)\varphi_{3,5}(b)+\varphi_{3,3}(a)\varphi_{3,4}(b)\\
&&+\varphi_{3,4}(a)\varphi_{3,3}(b)+\varphi_{3,5}(a)\varphi_{3,2}(b)\\
&&+\varphi_{3,6}(a)\varphi_{3,1}(b)+\varphi_{3,7}(a)\varphi_{3,0}(b)\\
&=&\sigma^3(a)d^3(b)+\sigma^2d(a)d^2\sigma(b)\\
&&+\sigma d\sigma(a)d\sigma d(b)+\sigma
d^2(a)d\sigma^2(b)\\
&&+d\sigma^2(a)\sigma d^2(b)+d\sigma d(a)\sigma d\sigma(b)\\
&&+d^2\sigma(a)\sigma^2 d(b)+d^3(a)\sigma^3(b).
\end{eqnarray*}
\end{example}

\begin{corollary} Let $d\sigma=\sigma d=d$. Then for each $a,b\in{\mathcal
D}$ we have
\[d^n(ab)=\sum_{r=0}^n {n\choose r} d^r(a)d^{n-r}(b).\]
\end{corollary}
\begin{proof} If $k$ represented in base $2$ has $r$ 1's,
then $\varphi_{n,k}=d^r$. But we have $n\choose r$ terms in the
summand with exactly $r$ 1's in representation of $k$.
\end{proof}

Note that by putting $\sigma=\iota$, we get the known results concerning ordinary derivations.

Our next result generalizes Theorem 3.2. As before, let $k$
is represented as $(c_n\ldots c_1)_2$ in base $2$. If the number
of $1$'s in this representation is $r_k$, we can construct
$2^{r_k}$ numbers $t$ with the property that $1$ occurs in $t$
only if the corresponding position at the representation of $k$
is $1$. More precisely, we can write
\[ T_k=\{t=(d_n\ldots d_1)_2~:\;\;\;\;d_i=1 \mbox{ implies } c_i=1 \mbox{ for each } 1\leq
i\leq n\}.\] To illustrate $T_k$'s, let $k=19=(10011)_2$. Then
\begin{eqnarray*}T_{19}&=&\{(00000)_2,(00001)_2,(00010)_2,(00011)_2,\\
&&(10000)_2,(10001)_2,(10010)_2,(10011)_2\}\\
&=&\{0,1,2,3,16,17,18,19\}.\end{eqnarray*} Here $T_k$ has $2^3=8$
elements.
\begin{lemma} Suppose $n,k$ are two natural numbers. Then

(i) $T_0=\{0\}, T_{2^n}=\{0,2^n\}$ and
$T_{2^n-1}=\{0,1,2,\ldots,2^n-1\}$;

(ii)
$T_{k}=T_{k-2^n}\cup(2^n+T_{k-2^n})=T_{k-2^n}\cup\{2^n+t~:~t\in
T_{k-2^n}\}$, provided that $2^n\leq k\leq 2^{n+1}-1$.
\end{lemma}
\begin{proof} (i) It is clear.

(ii) Let $k=(c_n\ldots c_1)_2$. We have $c_n=1$, since $k\geq
2^n$. Let $(d_n\ldots d_1)_2\in T_k$. If $d_n=0$ then $(d_n\ldots
d_1)_2=(d_{n-1}\ldots d_1)_2\in T_{k-2^n}$, and if $d_n=1$ then
$(d_n\ldots d_1)_2=2^n+(d_{n-1}\ldots d_1)_2$, where
$(d_{n-1}\ldots d_1)_2\in T_{k-2^n}$.
\end{proof}

\begin{definition} A linear mapping $\sigma$ is called a semi-endomorphism
if \[\sigma(a\sigma(b))=\sigma(a)\sigma^2(b)\] and
\[\sigma(ad(b))=\sigma(a)\sigma(d(b))\] for all $a, b \in{\mathcal
D}$. Obviously any endomorphism is semi-endomorphism.
\end{definition}

\begin{theorem} Let $\sigma$ be an endomorphism. Then for
each $n,k\in{\mathbb N}$ with $0 \leq k\leq 2^n-1$ and
$a,b\in{\mathcal D}$ we have
\[ \varphi_{n,k}(ab)=\sum_{\ell\in
T_k}\varphi_{n,\ell}(a)\varphi_{n,k-\ell}(b)\hspace{2cm}(***)\]
\end{theorem}
\begin{proof} We use induction on $n$. For $n=1$, if $k=0$
then $(***)$ is clear and if $k=1$ then $T_1=\{0,1\}$ and
\[
\varphi_{1,1}(ab)=d(ab)=d(a)\sigma(b)+\sigma(a)d(b)=\varphi_{1,1}(a)\varphi_{1,0}(b)+
\varphi_{1,0}(a)\varphi_{1,1}(b).\] Now suppose that $(***)$ is true
for $n$. For $0\leq k=(c_{n+1}c_n\ldots c_1)_2\leq 2^{n+1}$, two
cases occur.

{\bf Case1:} $1\leq k<2^n$.

In this case, $c_{n+1}=0$ and
$\varphi_{n+1,k}=\sigma\varphi_{n,k}$. Whence
\begin{eqnarray*} \varphi_{n+1,k}(ab)&=&\sigma\varphi_{n,k}(ab)\\
&=&\sigma(\sum_{\ell\in T_k}\varphi_{n,\ell}(a)\varphi_{n,k-\ell}(b))\\
&=&\sum_{\ell\in T_k}\sigma\varphi_{n,\ell}(a)\sigma\varphi_{n,k-\ell}(b)\\
&=&\sum_{\ell\in
T_k}\varphi_{n+1,\ell}(a)\varphi_{n+1,k-\ell}(b).\end{eqnarray*}

{\bf Case2:} $2^n\leq k<2^{n+1}-1$.

In this case, $c_{n+1}=1$ and so
$\varphi_{n+1,k}=d\varphi_{n,k-2^n}$. Thus
\begin{eqnarray*}\varphi_{n+1,k}(ab)&=&d\varphi_{n,k-2^n}(ab)\\
&=&d(\sum_{\ell\in T_{k-2^n}}\varphi_{n,\ell}(a)\varphi_{n,k-2^n-\ell}(b))\\
&=&\sum_{\ell\in T_{k-2^n}}[d\varphi_{n,\ell}(a)\sigma\varphi_{n,k-2^n-\ell}(b)\\
&&+ \sigma\varphi_{n,\ell}(a)d\varphi_{n,k-2^n-\ell}(b)]\\
&=&\sum_{\ell\in
T_{k-2^n}}[\varphi_{n+1,2^n+\ell}(a)\varphi_{n+1,k-2^n-\ell}(b)
\\&&+\varphi_{n+1,\ell}(a)\varphi_{n+1,k-\ell}(b)]\\
&=&\sum_{2^n\leq m\in
T_k}\varphi_{n+1,m}(a)\varphi_{n+1,k-m}(b)\\&&+\sum_{2^n>m\in
T_k}\varphi_{n+1,m}(a)\varphi_{n+1,k-m}(b)\\
&=&\sum_{m\in
T_k}\varphi_{n+1,m}(a)\varphi_{n+1,k-m}(b).
\end{eqnarray*}
\end{proof}

\begin{remark} Putting $k=2^n-1$ in the above theorem we get Theorem 3.2.
\end{remark}

The following theorem with $\sigma=\iota$ is a generalization of
Wielandt--Wintner theorem (cf. Theorem 2.2.1 of \cite{SAK}):

\begin{theorem} Let $\sigma$ be a bounded endomorphism on the Banach algebra ${\mathcal A}$, $d$ be
a bounded $\sigma$-derivation such that $d\sigma=\sigma d=d$ and
$d^2(a)=0$. Then $d(a)$ is a quasinilpotent, i.e. $r(d(a))= 0$.
\end{theorem}
\begin{proof} Using induction on $n$ we can establish that $d^n(a^n)=n!d(a)^n$ holds for all positive integer
$n$. Indeed if $d^{n-1}(a^{n-1})=(n-1)!d(a)^{n-1}$ then we infer
from Corollary 3.4 that
\begin{eqnarray*}
d^n(a^n)=d^n(a^{n-1}a)&=&\sum_{r=0}^n{n\choose
r}d^r(a^{n-1})d^{n-r}(a)\\
&=&n d^{n-1}(a^{n-1})d(a)+d^n(a^{n-1})\\
&=&n (n-1)!d(a)^{n-1}d(a)+d(d^{n-1}(a^{n-1}))\\
&=&n!d(a)^n+d((n-1)!d(a)^{n-1})=n!d(a)^n+0.
\end{eqnarray*}

Hence $r(d(a))=\displaystyle{\lim_{n\to\infty}}\|d(a)^n\|^{1/n}
=\displaystyle{\lim_{n\to\infty}}\|d^n(a^n)\|/n!)^{\frac{1}{n}}\leq
\frac{\|d\|\|a\|}{(n!)^{1/n}}=0$.
\end{proof}

We are ready to extend the Wielandt--Wintner theorem which states
there are no two elements $a$ and $b$ in a Banach algebra such
that $ab-ba=1$ (see Corollary 2.2.2 of \cite{SAK}).

\begin{theorem} Let $\sigma$ be a bounded endomorphism on the Banach algebra ${\mathcal A}$.
Then there are no three elements $a,b,c\in{\mathcal A}$ satisfying the following conditions such
that $a\sigma(b)-\sigma(b)a=c$:

(i) $\sigma(a)\sigma^2(b)-\sigma^2(b)\sigma(a)=a\sigma(b)-\sigma(b)a$,

(ii) $(\sigma^2(b)-\sigma(b))a=a(\sigma^2(b)-\sigma(b))$

(iii) $a\sigma(c)-\sigma(c)a=0$

(iv) $c$ is not quasinilpotent.
\end{theorem}
\begin{proof} Use the previous theorem with the inner
$\sigma$-derivation $d_a(u)=a\sigma(u)-\sigma(u)a, u\in {\mathcal
A}$. In fact, the conditions implies that $d_a^2(b)=0$ and so
$d_a(b)=c$ would be quasinilpotent which is a contradiction.
\end{proof}


\begin{thebibliography}{99}

\bibitem{A-M-N} Gh. Abbaspour, M. S. Moslehian and A. Niknam, Dynamical systems
on Hilbert $C^*$-modules, \textit{Bull. Iranian Math. Soc.}
\textbf{31(1)} (2005), 25--35.

\bibitem{A-R} M. Ashraf and N. Rehman, On $(\sigma-\tau)$-derivations in prime rings,
\textit{Arch. Math. (BRNO)} \textbf{38} (2002), 259--264.

\bibitem{B-M} C. Baak and M. S. Moslehian, On the stability of $\theta$-derivations on $JB^*$-triples,
\textit{Bull. Braz. Math. Soc.} \textbf{37} (2006), no. 1, 1--13.

\bibitem{BAT} C. J. K. Batty, Derivations on compact spaces, \textit{Proc. London Math. Soc.} \textbf{(3) 42} (1981),
no. 2, 299--330.

\bibitem {DAL} H. G. Dales, \textit{Banach algebras and automatic continuity},
London Mathematical Society Monographs, 24, Clarendon Press, Oxford,
2000.

\bibitem{B-R1} O. Bratteli and D. R. Robinson,  \textit{Operator Algebras and
Quantum Statistical Mechanics}, Vol. I, Springer-Verlag, New York,
1987.

\bibitem{B-R2} O. Bratteli and D. R. Robinson,  \textit{Operator Algebras and
Quantum Statistical Mechanics}, Vol. II, Springer-Verlag, New York,
1997.

\bibitem{H-K} M. Hongan and H. Komatsu, $(\sigma,\tau)$-derivations with invertible values,
 \textit{Bull. Inst. Math. Acad. Sinica} \textbf{15(4)} (1987), 411–-415.

\bibitem{KOM} H. Komatsu, On inner $(\sigma,\tau)$-derivations,  \textit{Math. J. Okayama Univ.} \textbf{23(1)} (1981),
33–-36.

\bibitem{M-M} M. Mirzavaziri and M. S. Moslehian, Automatic continuity of $\sigma$-derivations in
$C^*$-algebras,  \textit{Proc. Amer. Math. Soc.} \textbf{134}
(2006), no. 11, 3319--3327.

\bibitem{MOS1} M. S. Moslehian, Hyers--Ulam--Rassias stability of generalized derivations,
\textit{Intern. J. Math. Math. Sci.} \textbf{2006} (2006), 93942,
1–-8.

\bibitem{MOS2} M. S. Moslehian, Approximate $(\sigma-\tau)$-contractibility, \textit{Nonlinear Funct. Anal. Appl.}
(in press).

\bibitem{SAK} S. Sakai, \textit{Operator Algebras in Dynamical Systems}, Cambridge Univ. Press, 1991.
\end{thebibliography}
\end{document}